\documentclass[12pt]{article}

\title{A SIMPLE PRACTICAL HIGHER DESCENT FOR LARGE HEIGHT RATIONAL POINTS
ON CERTAIN ELLIPTIC CURVES}
\author{Allan J. MacLeod\\Dept. of Mathematics and Statistics\\
University of Paisley\\High St.\\Paisley\\Scotland\\PA1 2BE\\
(e-mail: macl-ms0@paisley.ac.uk)}

\begin{document}

\maketitle

\begin{abstract}

We consider the problem of determining a rational point $(x,y)$, of infinite
order, on the elliptic curve $y^2=x^3+ax^2+bx$. The well-known method of 
4-descent has proven very effective, and we give explicit formulae for the 
various stages.
For points of large height, however, there can still be a significant search
to be done. We describe a very simple idea which essentially doubles the 
height range which can be dealt with in a reasonable time. Currently, the
largest point found has height 51.15 or 102.3 (depending on which normalisation
you use). The report is written in an algorithmic style and should be
understandable by the interested amateur.

\end{abstract}

\newpage

\section{Introduction}
Many Diophantine problems can be reduced to determining points on elliptic 
curves. A fascinating example is from Bremner et al \cite{bgn}, 
where finding possible representations of $n$ as

\begin{displaymath}
n = ( x + y + z) \left ( \frac{1}{x} + \frac{1}{y} + \frac{1}{z} \right )
\end{displaymath}
is shown to be equivalent to finding points of infinite order on the
elliptic curve

\begin{displaymath}
E_n: u^2 = v^3 + (n^2-6n-3) v^2 + 16n v
\end{displaymath}

The rank of $E_n$ can be estimated from the L-series of the curve, 
assuming the Birch \& Swinnerton-Dyer conjecture, but this just provides 
evidence of existence. For
an actual representation, we need to explicitly find a rational point. If the
curve has rank greater than 1, it is usually fairly simple to find a point
by a trivial search. For rank 1 curves, however, this is often not feasible.
The L-series computations can be extended to give an estimate for the height
of a rational point of the curve. If this is large, then a simple search will
take far too long.

Recently, Silverman  \cite{silv} produced a very nice and effective procedure 
to determine such points. Based on a preprint of this paper,
the author coded the method in UBASIC, and used it 
on several families of elliptic curves from Diophantine problems. 
This usage has shown that the method starts to become
time-consuming for heights greater than about 10 in the Silverman 
normalisation - 20 in the alternative (from now on all heights will be
in the Silverman normalisation).

If the curve has a rational point of order 2, we can use a 4-descent procedure
to calculate the point, which we have found to be effective up to heights 
of about 15. By effective we mean that a point can usually be found in a
matter of minutes on a 200MHz PC with a program written in UBASIC. If one is
willing to wait several hours, the height barrier can obviously be increased.
It should be noted that Silverman's method works just as well on curves with
no rational points of order 2, and so is more general.

\section{Four-descent}
The following set of formulae describe the algebra necessary to perform
a general four-descent procedure for the elliptic curve

%eqn 1
\begin{equation}
y^2 = x^3 + a x^2 + b x
\end{equation}
so the curve is assumed to have at least one point of order 2. This 
discussion is purely concerned with computing a rational point, and not with
determining the rank, and so is not as general as Cremona's {\bf mwrank}
program, described in \cite{crem}. 

We can assume without loss of generality that $x = du^2/v^2$ and $y=duw/v^3$,
with $d$ squarefree and $(d,v)=1$, giving

%eqn 2
\begin{equation}
d^2 w^2 = d^3 u^4 + d^2 a u^2 v^2 + d b v^4
\end{equation}
which implies that $d|b$, so that $b=de$. There are thus a finite set of 
possible values for d, which we can work out easily from the factors of
$b$, unless $b$ is very large and difficult to factor. Remember that $d$
can be negative.

Consider first the simpler quadratic
%eqn 3
\begin{equation}
h^2 = d f^2 + a f g + e g^2
\end{equation}
and look for a solution $(h_0,f_0,g_0)$ either by searching or by more
advanced methods, and assume $g_0 \ne 0$.

If we find a solution, we can parameterise as follows. Define $x=f/g$ and
$y=h/g$, so that the equation is 

\begin{displaymath}
y^2=d x^2+a x+e
\end{displaymath}
then the line through $(x_0,y_0)=(f_0/g_0,h_0/g_0)$, with gradient $m$,
meets the curve again where

%eqn 4
\begin{equation}
x=\frac{a + d x_0 + m (m x_0 - 2 y_0)}{m^2-d}
\end{equation}

Assuming $m=p/q$ is rational, we simplify this to

%eqn 5
\begin{equation}
\frac{a g_0 q^2 + d f_0 q^2 + p (f_0 p - 2 h_0 q)}{g_0 (p^2 - d q^2)}
\end{equation}

To solve our original problem we look for values of the parameters
giving $f/g=u^2/v^2$, which leads to the two quadratics

%eqn 6
\begin{equation}
k u^2 = g_0 p^2 - g_0 d q^2
\end{equation}
and

%eqn 7
\begin{equation}
k v^2 = f_0 p^2 - 2 h_0 p q + (a g_0 + d f_0)q^2
\end{equation}
with $k$ squarefree.

The possible values of $k$ are those which divide the resultant of the 
two quadratics, which means that $k$ must divide the determinant

\begin{center}
\begin{math}
\begin{array}{|cccc|}
g_0&0&-d g_0&0 \\
0&g_0&0&-d g_0 \\
f_0&-2 h_0&a g_0+d f_0&0 \\
0&f_0&-2 h_0&a g_0 + d f_0 
\end{array}
\end{math}
\end{center}
which reduces to $g_0^4(a^2-4b)$. Let $k_0$ be a squarefree divisor (again 
possibly negative).

We search the first of these quadratics to find a solution $(u_0,p_0,q_0)$.
If we find one, another simple line-quadratic intersection
analysis characterises the solutions to the first quadratic as

%eqn 8
\begin{equation}
p = 2 u_0 k_0 r s - p_0 (g_0 s^2 + k_0 r^2)
\end{equation}
and

%eqn 9
\begin{equation}
q = q_0 ( g_0 s^2 - k_0 r^2)
\end{equation}

Substitute these into $k_0 * (7)$, giving the quartic

%eqn 10
\begin{equation}
z^2 = z_1 r^4 + z_2 r^3 s + z_3 r^2 s^2 + z_4 r s^3 + z_5 s^4
\end{equation}
with

%eqn 11
\begin{equation}
z_1 = k_0^3 (a g_0 q_0^2 + d f_0 q_0^2 + f_0 p_0^2 - 2 h_0 p_0 q_0)
\end{equation}

%eqn12
\begin{equation}
z_2 = 4 u_0 k_0^3 (h_0 q_0 - f_0 p_0)
\end{equation}

%eqn 13
\begin{equation}
z_3 = 2 k_0^2 (f_0 (2 u_0^2 k_0 - d g_0 q_0^2 + g_0 p_0^2) -
 a g_0^2 q_0^2)
\end{equation}

%eqn 14
\begin{equation}
z_4 = -4 u_0 g_0 k_0^2 (f_0 p_0 + h_0 q_0)
\end{equation}

%eqn 15
\begin{equation}
z_5 = g_0^2 k_0 (a g_0 q_0^2 + d f_0 q_0^2 + f_0 p_0^2 + 2 h_0 p_0 q_0)
\end{equation}

This quartic can then be searched for possible solutions. Having found one,
the various transformations eventually lead back to a point on the curve.

\section{COMPUTING}
The above formulae form the basis of a simple code written by the author
in UBASIC. This system is fast, simple, free, and runs on a multitude of PC's,
old and new. The major advantage, however, is that large numbers can be dealt
with by one single statement at the start of the program - the rest of the
code is standard Basic. Constructing a UBASIC code also leads to a 
structure which can easily be translated into Fortran, C, C++, etc using the
many multiple-precision packages available in these languages.

The basic structure of the algorithm is

\begin{flushleft}
find divisors d of b

for s1 = s1a to s1b

\hspace{2em}test if $h^2=df^2+afg+eg^2$ soluble with $|f|+|g|=s1$

\hspace{2em}if $(d,h_0,f_0,g_0)$ a solution then 

\hspace{2em}find divisors k of $g_0(a^2-4b)$

\hspace{2em}for s2 =2 to s2b

\hspace{4em}test if $k u^2 = g_0 p^2 - g_0 d q^2$ with $|p|+|q|=s2$ 

\hspace{4em}if $(k_0,u_0,p_0,q_0)$ a solution, then

\hspace{4em}form quartic equation from equations (10) to (15)

\hspace{4em}test if quartic is soluble, and if it is

\hspace{4em}for s3=2 to s3b

\hspace{6em}test if quartic is square for $|r|+|s|=s3$

\hspace{6em}if it is, determine point on curve, and stop

\hspace{4em}next s3

\hspace{2em}next s2

next s1

\end{flushleft}

If the curve has a rational point of moderate height, the above search-based
procedure works very well. We have found the method to be effective for points
with height of up to about 15, if we
select s2b=99 and s3b=199. It is impossible to predict in advance the best
choice of s1a and s1b.

For larger heights, the quartic search needs a larger value of s3b and
can take a considerable time. In many
cases, the quartic is insoluble so searching is futile. Cremona describes
how to test this, based on ideas from Birch and Swinnerton-Dyer \cite{bsd},
and we have implemented this test as a precursor to
searching. This method is probably the most advanced mathematically
in the whole procedure. It is also a vital time-saving procedure since
the majority of quartics are not soluble.

Cremona also describes various methods for searching the quartic
which reduce the time taken. We have chosen NOT to implement these, partly
to keep the code reasonably understandable by amateurs, and partly because 
UBASIC has restrictions on the number of variables. These methods are 
only effective time-savers if s3b is large.

For investigating families of curves where $a$ and $b$ are functions of $n$,
we found the simple search for $(h_0,f_0,g_0)$ could fail to find any 
solutions in the specified range. In such a case we can adapt the search as
follows. If

\begin{displaymath}
h^2 = d f^2 +a f g + b g^2 / d
\end{displaymath}
then 

\begin{displaymath}
d(2h)^2 = ( 2df+ ag)^2 - (a^2-4b)g^2
\end{displaymath}
and if we factor $a^2-4b = \alpha \beta^2$ with $\alpha$ squarefree,

%eqn 16
\begin{equation}
dH^2 = F^2 - \alpha G^2
\end{equation}
with $F = 2df+ag , \; G=\beta g, \; H = 2h$. From a $(F_0,G_0,H_0)$ solution,
we can recover integer solutions
to the original equation as $f_0=\beta F_0 - a G_0 , \; g_0 = 2 d G_0, \;
h_0 = d \beta H_0$.

In equation (6), it is clear that $(\pm u_0,\pm p_0, \pm q_0)$ satisfy the
equation. We do not, however, need to consider all 8 possible combinations of
sign. From (11) to (15), the sign of $u_0$ is irrelevant if we allow 
negative values of r or s. Similarly, $-p_0$ and $-q_0$ lead to the same 
quartic as $p_0$ and $q_0$, and the other two possible combinations give
the same quartic. Thus there are only two essentially different quartics,
and it is easy to show the relationship between them implies that they are
both soluble or both insoluble. It is worthwhile to search both (if soluble),
as the coefficients are different so the fixed search range might give a 
solution for one but not the other.

At several points in our code, we need to factor numbers to find divisors.
This is done by an extremely simple-minded search procedure, without recourse
to any modern factorisation techniques. So far, this has not had a severe 
effect on the performance of the code, but one could easily devise an 
elliptic curve where the factorisation time was dominant. The current code
has the advantage that is is short.

Finally, it should be noted that the method could easily find a torsion point
of the curve, and not a point of infinite order. Since the latter is
usually wanted, the code checks that the point found is not a torsion
point from a list of x-values provided by the user.

\section{Further Descent}
As stated, the above method begins to become time-consuming at certain
heights. In such cases, several workers have resorted to impressive
further descents to determine rational points. A classic example is the
work of Bremner and Cassels on the curve $y^2=x^3+px$, see \cite{bc}.
These methods seem to be very problem-dependent and to involve manipulations in 
algebraic number fields.

The author was interested in a general method of higher descent which does not
involve anything more than rational arithmetic - especially for use by the
many amateurs interested in Diophantine equations. The following method is
based on a remarkably simple idea. The author cannot believe this has not been
thought of before, but can find no direct reference to the idea.

The problem in the 4-descent for large heights is the determination of
$(r,s)$ values giving a square quartic. Suppose

\begin{equation}
z_1r^4+z_2r^3s+z_3r^2s^2+z^4rs^3+z_5s^4=(u_1r^2+u_2rs+u_3s^2)(v_1r^2+v_2rs
+v_3s^2)
\end{equation}
with $u_1,u_2,u_3,v_1,v_2,v_3$ integers,
then we can consider the following two quadratics

\begin{equation}
u_1r^2+u_2rs+u_3s^2=k_1m^2
\end{equation}

\begin{equation}
v_1r^2+v_2rs+v_3s^2=k_1t^2
\end{equation}
with $k_1$ squarefree. As before $k_1$ divides the resultant, which is

\begin{equation}
u_1^2v_3^2-u_1u_2v_2v_3-2u_1u_3v_1v_3+u_1u_3v_2^2+
u_2^2v_1v_3-u_2u_3v_1v_2+u_3^2v_1^2
\end{equation}

Suppose we find a solution $(k_1,t_1,r_1,s_1)$ to the second quadratic,
then we can parameterise as

\begin{equation}
 r=i^2k_1r_1-2ijk_1t_1+j^2(r_1v_1+s_1v_2)
\end{equation}

\begin{equation}
s=s_1(i^2k_1-j^2v_1)
\end{equation}
which, when substituted into $k_1*(18)$, requires the following quartic
to be square

\begin{equation}
c_1 i^4 + c_2 i^3 j + c_3 i^2 j^2 + c_4 i j^3 + c_5 j^4
\end{equation}
with

\begin{equation}
c_1=k_1^3(r_1^2 u_1+r_1 s_1 u_2 + s_1^2 u_3)
\end{equation}

\begin{equation}
c_2=-2 k_1^3 t_1 ( 2 r_1 u_1 + s_1 u_2)
\end{equation}

\begin{equation}
c_3=k_1^2 (4 k_1 t_1^2 u_1 + 2 r_1^2 u_1 v_1 + 2 r_1 s_1 u_1 v_2
+ s_1^2 (u_2 v_2 - 2 u_3 v_1))
\end{equation}

\begin{equation}
c_4=- 2 k_1^2 t_1 (2 r_1 u_1 v_1 + s_1 ( 2 u_1 v_2 - u_2 v_1 ) ) 
\end{equation}

\begin{equation}
c_5=k_1 (r_1^2 u_1 v_1^2 + r_1 s_1 v_1 ( 2 u_1 v_2 - u_2 v_1 ) + 
s_1^2 (u_1 v_2^2 - v_1 (u_2 v_2 - u_3 v_1)))
\end{equation}

This quartic can be tested for solubility and then searched.

\section{Computing the 8-descent}
The basic structure of the algorithm is obviously very similar to
the 4-descent structure.

\bigskip

\begin{flushleft}
find divisors d of b

for s1 = s1a to s1b

\hspace{2em}test if $h^2=df^2+afg+eg^2$ soluble with $|f|+|g|=s1$

\hspace{2em}if $(d,h_0,f_0,g_0)$ a solution then 

\hspace{2em}find divisors k of $g_0(a^2-4b)$

\hspace{2em}for s2 = 2 to s2b

\hspace{4em}test if $k u^2 = g_0 p^2 - g_0 d q^2$ with $|p|+|q|=s2$ 

\hspace{4em}if $(k_0,u_0,p_0,q_0)$ a solution, then

\hspace{4em}form quartic equation from equations (10) to (15)

\hspace{4em}test if quartic soluble, and if so
 
\hspace{4em}test if quartic can be factored into 2 integer quadratics

\hspace{4em}if so, find possible values of $k_1$ 

\hspace{4em}for s3=2 to s3b

\hspace{6em}test if $k_1 t^2 = v_1 r^2 + v_2 rs + v_3 s^2$ with $|r|+|s| \le
s3$

\hspace{6em}if $(k_1,t_1,r_1,s_1)$ a solution, form quartic (20)

\hspace{6em}test if quartic soluble, and if it is

\hspace{6em}for s4=2 to s4b

\hspace{8em}test if $c_1i^4+c_2i^3j+c_3i^2j^2+c_4ij^3+c_5j^4$ square

\hspace{8em}if it is, use various transformations to find point and stop

\hspace{6em}next s4

\hspace{4em}next s3

\hspace{2em}next s2

next s1

\end{flushleft}

\medskip
The factorisation of the quartic is accomplished by considering the equations:

\begin{center}
\begin{math}
\begin{array}{ccl}
z_1&=&u_1 v_1 \\
z_5&=&u_3 v_3 \\
z_2&=&u_1 v_2 + v_1 u_2 \\
z_4&=&u_3 v_2 + v_3 u_2 \\
z_3&=&u_1 v_3 + u_2 v_2 + u_3 v_1
\end{array}
\end{math}
\end{center}

We thus have to factorise $z_1$ and $z_5$. For each possible splitting,
we solve the third and fourth equations for $u_2$ and $v_2$. If the
solutions are integers we test whether the 6 values satisfy the last equation.
We found that quadratics which themselves split into 2 linear factors lead
to torsion points, so we test whether the quadratics have rational roots,
and reject those which do.

\section{Numerical Examples}
In this section we give some examples of the use of the 8-descent method,
with some specimen timings on either a 200MHz or 300MHz PC.
It is clear that timings are dependant on the search parameters, so we specify
the values of (s2b,s3b,s4b). s1a is set to 2 and we search until
a solution is found.

\bigskip

(a) Congruent numbers N are integers which can be the area of a rational
right-angled triangle. Finding the sides of such a triangle is equivalent
to finding a point of infinite order on the elliptic curve
$y^2=x^3-N^2x$. 

One of the most famous of these numbers is $N=157$, since
it forms the basis of an impressive diagram in Chapter 1 of Koblitz's
book \cite{kob}, where the sides involve numbers with 20-30 digits. An L-series
calculation gives the height of a point as $27.3$. Running the program
with s2b=s3b=99 and s4b=199 finds the point with x-coord

\begin{displaymath}
\frac{-1661 \, 3623 \, 1668 \, 1852 \, 6754 \, 0804}{28 \, 2563 \,
0694 \, 2511 \, 4585 \, 8025}
\end{displaymath}
in 18.8 secs (200MHz). With s4b=99 it takes only 7.0 secs, while for s4b=499
it takes 31.5 secs. But, with s4b=599 the time goes down to 28.3 secs,
and if s4b=699 the time is only 9.0 secs. The variation is due to the 
number of quartics which need to be searched. Up to 499 it takes 22
quartics, but 599 needs only 4, while 699 finds a point on the first.
The important point is that it takes less than 1 minute to find such
a large point. Even an ancient 80387 machine finds this point in less 
than 5 minutes.

The author has used this technique in a search for actual solutions to
the congruent number problem for $1 \le N \le 1999$.
The current method has led to the completion of a table of
solutions for $[1,499]$. The largest height encountered was $51.15$ for 
$N=367$. Searching with s2b=s3b=99 and s4b=9999, the following
solution was found on a 300MHz machine in 17595 secs.

\begin{displaymath}
x = \frac{-367 (4 \, 9695 \, 3629 \, 6085 \, 1360 \, 8777)^2}
{(163 \, 8216 \, 8821 \, 6485 \, 0643 \, 1464)^2}
\end{displaymath}

It is perfectly possible that the modular form Heegner-point method of
Elkies \cite{elk} could find this point much faster, but this method is 
much more difficult for the non-expert to understand. Elkies method also does
not generalise to other families of curves.

Since $x^3-N^2x=x(x-N)(x+N)$, we can change the origin to produce the
two equivalent elliptic curves $y^2=x^3 \pm 3Nx^2 + 2N^2x$. They
may be equivalent mathematically, but their performance computationally
is not the same. The author has found several solutions to the congruent
number problem from these curves having been unsuccessful with the
original curve. This also happens in other families of elliptic
curves with 3 rational points of order 2.

\bigskip

(b) The paper of Bremner et al mentioned in the introduction has a
representation of $N=564$ in its abstract. The L-series computation
predicts a point with height $38.01$. A 300MHz PC found a solution
in 4497 secs. with s2b=s3b=99 and s4b=9999. This solution
leads to a representation with

\medskip

x = 32736 87951 95203 44322 22320 98479 60433 77911 47254 01060

y = 53 58267 18225 66098 96868 10234 90522 46809 90105 26717

z = - 1158 25525 22781 02629 66659 36639 59067 36616 11576 01937

\bigskip
(c) The paper of Bremner and Cassels describes the finding of a 
point on $y^2=x^3+877x$. The L-series predicts a height of $48.0$.
A 300MHz PC finds the following x-coordinate in 51874 seconds with
s2b=s3b=99 and s4b=12999.
\begin{displaymath}
x = \frac{877 \, ( 7884 \, 1535 \, 8606 \, 8390 \, 0210 )^2}
{(6 \, 1277 \, 6083 \, 1879 \, 4736 \, 8101 )^2}
\end{displaymath}

\bigskip
(4) The final example is included for historical reasons, since it was
the first time that the author tried the idea of factorising the
quartic. Since this was the start of the current research, clearly
the computer programs used in the previous examples were not available.
The computations were done by simple searches and algebraic manipulation
with Derive.

The problem comes from the diophantine problem of finding an
integer triangle with base/altitude = n. For $n=79$, we consider the
equation $y^2=x^3+6243x^2+x$. The L-series computations suggested rank 1,
but with a height of over 40 for a rational point. The 2-isogenous curve
is $y^2=x^3-12486x^2+38975045x$, which was indicated to have a point
with height about 20.

The author selected to try $d=5$, which meant looking first for
solutions of $h^2=5f^2-12486fg+7795009g^2$. A very simple search program
quickly finds $f=93, g=1, h=2584$, which means (6) and (7), are
\begin{displaymath}
k u^2 = p^2 - 5 q^2
\end{displaymath}
\begin{displaymath}
k v^2 = 93 p^2 - 5168 pq -12021 q^2
\end{displaymath}
where $k=\pm 1$. We selected $k=1$.

It is possible to parameterise the first as $p=r^2+5s^2$,$q=2rs$, which gives
\begin{displaymath}
v^2=93r^4-10336r^3s-47154r^2s^2-51680rs^3+2325s^4
\end{displaymath}
which Derive fairly easily factors as
\begin{displaymath}
v^2=(3r^2-340rs-775s^2)(31r^2+68rs-3s^2)
\end{displaymath}

The two quadratic factors form the basis of (18) and (19), and we find
that $k_1|158$. Picking $k_1=1$, we found several $(r,s)$ solutions
to $3r^2-340rs-775s^2=t^2$, which lead to parameterisations for
r and s, but which lead to insoluble quartics. We then tried
$k_1=158$, and found the solution $r=9,s=-1,t=4$, which gives the
parameterisation
\begin{displaymath}
\frac{r}{s} = \frac{-(1422i^2+1264ij+367j^2)}{158i^2-3j^2}
\end{displaymath}
and hence to the quartic
\begin{displaymath}
316^2(74892 i^4 +154840 i^3j+123789i^2j^2+45916ij^3+6725j^4)
\end{displaymath}
which has to be square.

Thus quartic was everwhere soluble, so a search was started which quickly
found the solution $i=151, j=-158$, which lead back to a point on the
2-isogenous curve with
\begin{displaymath}
x=\frac{2836 8499 3467 6319 5139 0020}{4689 8490 9449 9234 0041}
\end{displaymath}
and thus to a point on the original curve with
\begin{displaymath}
x=\frac{2654 7926 1289 1944 1996 8505 1867 1143 3025}
{1705 4187 5947 7256 7676 9862 5643 5806 2336}
\end{displaymath}

For interested readers, this point leads to the triangle with sides

\begin{math}
\begin{array}{crr}
\, & \, & \, \\
A =& 1465 86997 18477 82318 35321& 97194 40069 87886 57474 85658 \\ 
\, &           \,                & 64108 26213 28674 16311 64960 \\
\, & \, & \, \\
B =& 892 76765 34887 48588 76033 & 62942 70957 75073 78277 30811 \\
\, &           \,                & 86659 99941 08625 53894 71249 \\
\, & \, & \, \\
C =& 573 59536 91823 05619 55378 & 66267 79319 29261 59738 76797 \\
\, &           \,                & 12797 54707 31247 71171 08209
\end{array}
\end{math}

\section{Further Work}
The reasonable level of success of this method for several families of
elliptic curves provided the impetus to write this report, but there is still
much work to be done.

\begin{enumerate}
\item Translate the UBASIC code into a compiled high-level language so
that the method can be run on non PC's, especially UNIX workstations.

\item Try further descents on the last quartic produced. The author has
experimented with this idea, but initial results are disappointing in
that the method seems to be finding multiples of a generator rather than
the generator itself. Since the multiples have much larger height it is
currently no benefit to try these extra descents.

\item Try to underpin the method with some theory. The method was developed
by trying an idea, seeing it work, and improving it. The method does not 
always work, but does provide another tool in the investigator's toolbox.
\end{enumerate}

\end{document}